\input amstex
\documentstyle{amsppt}
\magnification1100

%\hsize=5.3truein
%\hoffset=-.5truein
%\vsize=9.7truein
%\voffset=-.5truein
\NoBlackBoxes
\catcode`\@=11
\redefine\logo@{}
\catcode`\@=13
\font\das=cmss10
 %scaled\magstep1

\font\idas=cmssi10
\font\bdas=cmssbx10

\addto\tenpoint{\advance\baselineskip by 2pt}\tenpoint

\def\ct#1\endct{\title#1\endtitle}

\def\at#1\endat{\title#1\endtitle}

\def\au#1\endau{\author#1\endauthor}

\def\af#1\endaf{\address#1\endaddress}
\def\da#1\endda{\date#1\enddate}
\def\dd#1\enddd{\dedicatory#1\enddedicatory}
\def\su#1\endsu{\subjclass#1\endsubjclass}
\def\kw#1\endkw{\keywords#1\endkeywords}
\def\tk#1\endtk{\thanks#1\endthanks}
\def\ab#1\endab{\abstract#1\endabstract}
\begingroup 
\let\head\relax \let\subhead\relax \let\subsubhead\relax
\outer\gdef\ah#1\endah{\head#1\endhead}
\outer\gdef\bh#1\endbh{\subhead#1\endsubhead}
\outer\gdef\ch#1\endch{\subsubhead#1\endsubsubhead}
\endgroup

\topmatter
\title
\das On just infinite pro-$p$-groups and\\arithmetically profinite extensions of local fields
\endtitle
\author
\das Ivan  Fesenko 
\endauthor
\rightheadtext{{\das Just infinite pro-$p$-groups and APF extensions}}
\leftheadtext{{\das I. Fesenko}}
\endtopmatter

\document
\predefine\til{\~}
\def\~#1{\relax\ifmmode\widetilde{#1}\else\til{#1}\fi}

\redefine \le{\leqslant}
\redefine \ge{\geqslant}
\define \wt#1{\mathaccent"0365{#1}}
\define \wh#1{\mathaccent"0362{#1}}

\define \Aut{\operatorname{Aut}}

\medskip

\das

Let $F$ be a local field of characteristic $p>0$ with perfect residue field
$k$. The wild group $R=\Aut_1 F$ is the group
of wild  continuous automorphisms 
$\{\sigma: (\sigma-1)\Cal O_F\subset 
\Cal M_F^2\}$ of the local field $F$. A choice of a prime element $t$ 
of a local 
field $F$ determines an isomorphism of  $\Aut_1 F$ and
the group  of formal power series $f(t)=t+a_2t^2+\dots$ with coefficients from $k$ with respect to the composition
$(f\circ g)(t)=f(g(t))$. We shall write $R_k$ for $R$ to specify the residue field.

The group $R_k$ plays an important role in the theory of pro-$p$-groups
being one of the main known representatives of the class of
just infinite pro-$p$-groups of finite width which are not
$p$-adic analytic, see for instance [18]. 
The group $R_k$ is known under various names including 
the Nottingham group and  the Richard group. 

It has been investigated by group theoretical methods (D. Johnson, I. York, A. Weiss, C.~Lee\-d\-ham-Green,
A. Weiss, A. Shalev, R. Camina, Y. Barnea, B. Klopsh) and number theoretical methods
(Sh. Sen, J.-M. Fontaine, J.-P. Wintenberger, F. Laubie, J. Lubin, H.-C.~Li).

In this paper we apply Fontaine--Wintenberger's  theory of fields of norms 
to study the structure of the wild group $R_k$.
In particular we  provide a new short proof of R.~Camina's theorem which says
that every countably based pro-$p$-group
(i.e. with countably many open sugroups) is isomorphic to
a closed subgroup of $R_{\,\Bbb F_p}$.
In sect. 6 we deal with some specific subgroups $T$ of $R_{\,\Bbb F_p}$;
they have  remarkable properties, in particular the commutator subgroup
is unusually small.  Realizing the group $T$ 
as the Galois group of arithmetically
profinite extensions of $\goth p$-adic fields
we answer affirmatively in sect. 7  Coates--Greenberg's problem
on deeply ramified extensions of local fields.
In final sect. 8 we show that the wild group is not analytic over
 commutative complete local noetherian
integral domains with finite residue field of characteristic $p$.

\bigskip

\subhead 1  \endsubhead Let $v_F$ be the discrete valuations of $F=k((t))$.
We use simultaneously two interpretations of $R=R_k$:
as formal power series $t+a_2t^2+\dots$ with respect to the composition
and as wild automorphisms of $F$.
When we use formal power series $f,g$ their multiplication is denoted by $f\circ g$,
when we use automorphisms $\sigma,\tau$ their product in $R$ is denoted by $\sigma\tau$.

For a formal power series $f(t)$ denote
$i(f(t))=\min \{i\ge 2: a_i\not=0\}-1$. 
Let
$R_i=\{f(t):i(f(t))\ge i\}$.

Denote $[\sigma,\tau]=\sigma\tau\sigma^{-1}\tau^{-1}$.
The following property of commutators will be often in use
$$[\sigma\tau,\rho]=[\sigma,[\tau,\rho]][\tau,\rho][\sigma,\rho].\tag"\das($\circ$)"$$

As usual, we  denote by dots terms of higher order.
For $a,b\in k$ 
$$[t+at^i,t+bt^j]=t+ab(i-j)t^{i+j-1}\,+\dots.$$ 
Therefore,  $R_i$ are normal subgroups of $R$.
The group $R$ being the projective limit of finite $p$-groups
$R/R_i$ is a pro-$p$-group.
From the commutator formula one immediately deduces that
 $[R_n,R_m]=R_{n+m}$ if $m-n$ is relatively prime to $p$
and  $[R_n,R_m]=R_{n+m+1}$ if $m-n$ is divisible by $p$, see for
example [7, Prop. 12]. In addition, $R_n^p\le R_{np}$, see the 
first part of the  proof of the Proposition in sect. 6.
Thus, $[R,R]=[R,R]R^p=R_3$. 
The group $R_{\,\Bbb F_p}$ is a pro-$p$-group with 2 generators,
generated by any two elements of $R_2\setminus R_3$ and $R_3\setminus R_4$,
hence by $t+t^2$ (of infinite order) and $t/(1-t)$ (of order $p$).

Moreover,
the group $R$ is a so called hereditarily just infinite group:
every nontrivial normal closed subgroup $G$ of an open subgroup is open.
Indeed, by the commutator formula the set $H=[t+at^{i+1}+\dots, R_j]$ has the  property
$R_u\le R_{u+1}H$ for $u\ge j+i, (p,u-i)=1$. Then for an odd $p$
and sufficiently large $l$ the group $G$ contains some $t+at^l+\dots$ and $t+at^{l+1}+\dots$, so $G$ contains $R_w$ for sufficiently large $w$.
For $p=2$ use in addition the property $(t+at^i)\circ(t+at^i)=t+a^2t^{2i-1}+\dots$.

\medskip

\subhead 2 \endsubhead For a closed subgroup $G$ of $R$ put
$G_i=G\cap R_i$ for $i\in\Bbb N$, $G_x=G_{\lceil x\rceil}$
for $x\in\Bbb R$. Denote
$\varphi_G(x)=\int_0^x\frac{d\,y}{|G:G_y|}$.
The group $G$ is called an arithmetically profinite
subgroup of $R$ if $\lim_{x\to+\infty}\varphi_G(x)=+\infty$, see
[25]. If this is the case, define $\psi_G(x)$ as the inverse function to
$\varphi_G(x)$ and put $G(x)=G_{\psi_G(x)}$.
The points of discontinuity of the derivative of $\varphi_G$ are called breaks of $G$; 
the points of discontinuity of the derivative of
$\psi_G$  are called upper breaks of $G$.

A theorem of Sen [23] says that the  subgroup of $G$
generated by an element $\sigma$ of infinite order is
arithmetically profinite and $i(\sigma^{p^n})\equiv i(\sigma^{p^{n-1}})
\mod p^n$. For generalizations and other proofs see  [27], [20], [19].

\medskip

\subhead 3 \endsubhead Let $K$ be a local field with perfect residue field $k=k_K$ of characteristic $p$. A Galois %totally ramified $p$-
extension $L/K$ is called arithmetically profinite if the upper ramification groups $G(L/K)^x$ are open in $G(L/K)$ for every $x$.
Equivalently, $L/K$ is  arithmetically profinite if
it has finite residue field extension and 
the Hasse--Herbrand function $h_{L/K}(x)=\lim h_{E/K}(x)$ takes real values
for all real $x\ge0$ where $E/K$ runs through all finite subextensions in $L/K$,
see [28], [9, Ch. III, sect. 5].

For an infinite Galois arithmetically profinite extension $L/K$ the field of norms $N=N(L/K)$ is the set of all norm-compatible sequences
$$\{(a_E):a_E\in E^*, E/K \text{\das \ \ is a finite subextension of $L/K$}\}$$
and zero, such that the multiplication is componentwise and the addition
$(a_E)+(b_E)=(c_E)$ is defined as $c_E=\lim_{M} N_{M/E}(a_M+b_M)$
where $M$ runs through all finite subextension of $E$ in $L$.
For the properties of
the fields of norms see [25], [28], [9, Ch. III, sect. 5].

In this paper with the last paragraph being excluded all Galois arithmetically profinite extensions are totally ramified $p$-extensions;
therefore the Galois group consists of wild automorphisms only.

The field $N$ is a local field of characteristic $p$ with the residue field $k_L$ and a prime element $t=(\pi_E)$ which is a sequence of norm-compatible
prime elements of finite subextensions of $L/K$.
Every automorphism $\tau$ of $L$ over $K$ being wild induces a wild automorphism $\sigma$ of the field of norms: $\sigma ((\pi_E))=(\tau \pi_E)$.

A Galois infinite subextension of a Galois arithmetically profinite extension
is arithmetically profinite.
Let $F/K$ be a Galois totally ramified $p$-extension and
$F$ contain $L$ which is an arithmetically profinite extension of $K$.
If $|F:L|<+\infty$, then $F/K$ is an arithmetically profinite extension.
%One can find a finite extension $M$ of $K$ such that $ML=F$.
%Then $N_{ME_1/ME_2}(a)=N_{E_1/E_2}(a)$ for $a\in E_1$ and sufficienlty large %$E_2$. Hence 
The field of norms $N(L/K)$ can be identified with a subfield of $N(F/K)$; 
%via $(a_E)\to (b_{E'})$ with $b_{ME}=a_E$ for sufficiently large $E$
%and  $b_{E'}=N_{ME/E'}b_{ME}$ for $E'\subset ME$. T
the extension
$N(F/K)/N(L/K)$ is an extension of local fields. 
If $F$ is a Galois extension of $L$, then
one defines $N(F,L/K)$ as the compositum of all $N(F'/K)$ where
$F'$ runs through Galois extensions of $K$ in $F$ with $|F':L|<+\infty$.
One of the central theorems of the theory of fields of norms says
that the absolute Galois group of $N(L/K)$ coincides with 
$G(N(L^{\text{\rm sep}},L/K)/N(L/K))$ and the latter is isomorphic to
$G(L^{\text{\rm sep}}/L)$, see [28, 3.2.2].

The functor of field of norms $W=W_{L/K}$ associates to an infinite Galois arithmetically profinite extension $L/K$ its field of norms $N(L/K)$ and  the closed
arithmetically profinite subgroup $G$ of the group $R_k=\Aut_1 N(L/K)$ which is the image of the Galois group of $L/K$;
the upper ramification filtration $G(L/K)^x$ is mapped onto the filtration
$G(x)$ of $G$, see [28, 3.3].

For a finite Galois totally ramified $p$-extension $N/K$ of local fields
of characteristic $p$ the Galois group $G(N/K)$ is isomorphic to a 
subgroup of $R_k=\Aut_1N$.
The extension $F/K$ is arithmetically profinite if and only if
$N(F,L/K)$ is an arithmetically profinite extension of
$N(L/K)$; then the image of $G(N(F,L/K)/N(L/K))$ under
$W_{N(F,L/K)/N(L/K)}$ in $R_k$ coincides with the image of
$G(F/L)$ as a closed subgroup in $G(F/K)$ in $R_k$ under
$W_{F/K}$, see [28, 3.2].

\medskip

\subhead 4 \endsubhead A theorem of Wintenberger ([27]) says that for every
abelian closed subgroup
$G$ of the group $R_k$ there exists a Galois arithmetically profinite extension
$L/K$ of local fields such that $W(L/K)=(k((t)),G)$.

For example, the group topologically generated by an element $\sigma$ of infinite order in $R$ comes from an arithmetically profinite  $\Bbb Z_p$-extension $L/K$. It is easy to deduce that the sequence $i(\sigma^{p^n})/p^n$ is increasing.
Denote 
$pe/(p-1)=\lim i(\sigma^{p^n})/p^n$.
Then either $e=+\infty$ or $e\in\Bbb N$. In the first case $K$ is of positive characteristic,
in the second case $K$ is of characteristic 0 and its absolute ramification index is $e$. By the proposition in sect. 6 
$e(\tau)=pi$ is finite for $\tau(t)=t+t^{1+pi}$.

An observation due to Fontaine [11] is that $e=+\infty$ if and only if
$\sigma$ belongs to the topological closure of the torsion 
(the set of torsion elements) of $R$.
Indeed, assume that the group $G$ topologically generated by
$\sigma$ comes from a Galois arithmetically profinite extension $L/K$ of fields of characteristic $p$ with a generator $\tau$.
Let $K_n$ be the subextension of $L$ of degree $p^n$ over $K$.
Map $K_n$ isomorphically onto $N$ by sending a prime element $\pi_{K_n}$ of a norm-compatible sequence of prime elements of finite subextensions
in $L/K$ to $t$. Let $N_n\subset N$ be the image of $K$ under this homomorphism,
and let $\sigma_n$ (of order $p^n$) be the image of $\tau$. Then $i(\sigma\sigma_n^{-1})$ tends to $+\infty$ when $n$ grows. 
Conversely, if $\sigma$ is the limit of a sequence of automorphisms
$\sigma_n$ of finite order, then the upper breaks $u_i$ of $G=\overline{(\sigma)}$
satisfy $u_{i+1}\ge u_i^p$ (see for instance [17]), therefore $e=+\infty$. 

Wintenberger and Laubie studied $p$-adic Lie subgroups in $R$ which are in the image of the functor $W$,  see [25], [26], [16].

\medskip

\subhead 5\endsubhead The wild group $R$ is not $p$-adic analytic,  since
for instance 
for every $n$ relatively prime to $p$ there is
$\sigma\in R_{n}\setminus R_{n+1}$ such that $\sigma^p=e$ 
(it suffices to observe that given a natural number relatively prime to $p$ there a cyclic totally ramified extension of degree $p$
of a local field of characteristic $p$ 
with the ramification break equal to that number).
Another way to argue is to use the property of $p$-adic analytic groups
to contain an open subgroup of finite rank (i.e. an open subgroup
for which the supremum of the number of generators of its closed subgroups 
is not infinity), see [6, Cor. 9.36]). The group
$R$ doesn't contain an open subgroup of finite rank, since the 
number of generators of $R_i$ tends to infinity when $i$ tends to infinity.

For more properties of $R$ see Remarks 1, 2 in sect. 6. 

\proclaim{\bdas Proposition} \idas Let $G$ be a countably based pro-$p$-group.
Then there exists a Galois arithmetically profinite extension $L$ of \ 
$\Bbb F_p((X))$ such that $G(L/\Bbb F_p((X)))$ is isomorphic to $G$.
 \endproclaim \das
\demo {\idas Proof}\das   
Let $G=\varprojlim G_i$ where $G_i$ are finite pro-$p$-groups
and $|G_{i+1}:G_i|=p$. Assume that $G_i$ is isomorphic to
$G(K_i/\Bbb F_p((X)))$. It is well known that the pro$p$-part of the absolute
Galois group of $\Bbb F_p((X))$ is a free countably generated pro-$p$-group.
Hence the imbedding problem 
$$(G_{i+1}\to G_i=G(K_i/\Bbb F_p((X))))$$ 
has a solution $K_i(\beta)$ with $\beta^p-\beta=\alpha\in
K_i$,  see  for instance [12, Th.1'].
Following the method of Camina [3] replace $\alpha$ by $\alpha_1=\alpha+c$ with $c\in \Bbb F_p((X))$.
Then $K_{i+1}=K_i(\beta_1)$ with $(\beta_1)^p-\beta_1=\alpha_1$
is a solution of the same imbedding problem and 
the ramification break of $K_{i+1}/K_i$ can be made arbitrarily large by 
choosing $c$ with $v_{\Bbb F_p((X))}(c)$ being sufficiently negative 
and relatively prime to $p$.
Therefore one can construct an arithmetically profinite extension
$L/\Bbb F_p((X))$ as desired.
\qed\enddemo 
\das

\proclaim{\bdas Corollary 1 (Camina)} \idas Every countably based pro-$p$-group
is isomorphic to a closed subgroup of $R_{\,\Bbb F_p}$.
\endproclaim \das
\demo {\idas Proof}\das
Apply the functor $W$ to the extension  $L/\Bbb F_p((X))$.
\qed\enddemo 
\das

\proclaim{\bdas Remark 1} \das According to the proof given in this paper
every countably based pro-$p$-group
is isomorphic to infinitely many different closed arithmetically profinite subgroups in the closure of the torsion of $R_{\,\Bbb F_p}$.
Note that if $\tau\in  G(L/\Bbb F_p((X)))$ is of infinite order, then
the fixed field $L_{\tau}$ of $\tau$ is an arithmetically profinite
extension of $\Bbb F_p((X))$ and the image of $\tau$ in $R_{\,\Bbb F_p}$
can be identified with the image of $\tau\in G(N(L, L_{\tau}/\Bbb F_p((X)))/
N(L_{\tau}/\Bbb F_p((X))))$ in $R_{\,\Bbb F_p}$. The latter belongs to the closure of the torsion of $R_{\,\Bbb F_p}$ as was indicated in the previous section.
Varying the set of upper ramification breaks as in the proposition  
every infinite countably based pro-$p$-group  can be  embedded
in infinitely many ways into $R$; the images have different sets of breaks.
\endproclaim

\proclaim{\bdas Remark 2} \das In Camina's proof every finitely generated pro-$p$-group is realized as the Galois group
of a totally ramified $p$-extension with specific properties of
its ramification breaks, then it is embedded into
$R_{\,\Bbb F_p}$ (actually in the closure of the torsion).
Then Lubotzky--Wilson's theorem is applied {\das(}see Corollary 2{\das)} to handle the general case of a  countably based pro-$p$-group.
One can use  Example 2.4 of [8] to show that the closed subgroups
of $R$ given by Camina's construction are not in general arithmetically profinite subgroups of $R$. 
\endproclaim

\proclaim{\bdas Remark 3} \das In discussions with D. Segal and B. Klopsch
we have observed that every closed subgroup $G$ of $R$ which is in the image of the functor of fields of norms has Hausdorff dimension (for the definition see [1]) equal to zero. Indeed, the nondecreasing sequence of the set of breaks $(s_i)$ of $G$
satisfies $\sum (s_i-s_{i_1})/p^i=+\infty$, hence $\lim\inf i/s_i=0$.
Closed subgroups of $R$
produced in Camina's construction have Hausdorff dimension zero as well.
\endproclaim

\das

The closure of the  torsion of $R$ is different from $R$, since
every automorphism with finite $e$ (see the previous section)
doesn't belong to the closure of the torsion of $R$
(by the way, the closure of the group generated by the torsion of $R$ coincides
with $R$ for $p>2$). 
The same arguments as in the proof
of Corollary 1 show that every closed non-pro-$p$-cyclic
subgroup of $R$ which is in the image of $W$ is inside the closure
of the torsion of $R$. 
Hence there is an infinite chain of closed subgroups of $R$:
$G_1=R>G_2>\dots$ such that  all $G_i$
are isomorphic to each other and
each next is contained in the closure of
the torsion of the previous one.

\proclaim{\bdas Corollary 2 (Lubotzky--Wilson)} \idas There is a pro-$p$-group with 2 generators which contains as a closed subgroup every countably based pro-$p$-group.
\endproclaim \das
\demo {\idas Proof}\das
The group $R_{\,\Bbb F_p}$ does.
\qed\enddemo 
\das

\proclaim{\bdas Problem} \idas Given a free pro-$p$-group $G$ with finite number of generators
does there exist a Galois arithmetically profinite extension $L$ of
$K$, $|K:\Bbb Q_p|<+\infty$, such that $G(L/K)$ is isomorphic to $G$?
\endproclaim \das

The affirmative answer will imply that for every finitely generated
pro-$p$-group $G$ there is a closed subgroup inside the group
$R_{\Bbb F_q}$
isomorphic to $G$ which comes via the functor of fields of norms
from a Galois arithmetically profinite extension of local number fields.

We will show in the next section that this is true for specific closed subgroups
$T$ of $R$ which are different from pro-$p$-cyclic groups.

\medskip
\subhead 6 \endsubhead For $m\ge2$ define the following closed 
subgroups in the wild group $R=R_{\,\Bbb F_p}$
$$S_m=\biggl\{\sum_{i\ge0}a_it^{1+mi}:a_0=1,a_i\in\Bbb F_p\biggr
\}.$$
For $m$ relatively prime to $p$ the group of principal units
$1+t^m\,\Bbb F_p[[t^m]]$ is uniquely $m$-divisible,
therefore one can associate to an element $\sigma\in R$
considered as a wild automorphism of $\Bbb F_p((t^m))$:
$\sigma(t^m)=t^m f(t^m)$ with $f(t)\in 1+t\Bbb F_p[[t]]$
an automorphism $\tau\in S_m$ considered as a wild automorphism of $\Bbb F_p((t))$: $\tau(t)=t\root{m}\of{f(t^m)}$. 
Hence $S_m$ is isomorphic to $R$
and $S_{mp^r}$ is isomorphic to $S_{p^r}$.

For $r\ge 1$ denote $q=p^r$, $T=S_{q}$ and put $T_i=\{f(t)\in T:
f(t)\in t+t^{1+qi}\,\Bbb F_p[[t]]\}$.

\proclaim{\bdas Proposition}\idas  

{\das (1)} If $\sigma\in T_i\setminus T_{i+1}$ then $\sigma^p\in
T_{pi}\setminus T_{pi+1}$; the intersection of $T$ with the
closure of the  torsion of $R$ is trivial.

{\das (2)} $[T_i,T_i]\le T_{(q+1)i+1}$, the group 
$T_i/T_{(q+1)i}$ is abelian of exponent $pq$.
A nontrivial normal closed sugroup of an open
subgroup of $T$ is open {\das (}$T$ is hereditarily just infinite{\das)}. 

{\das (3)} $[T,T]T^p>T_{q+2}$ and $T$ has not more than 
$q$ generators. 

{\das (4)} $T$ is not $p$-adic analytic.

\endproclaim \das
\demo {\idas Proof}\das

(1) For $\alpha\in F=\Bbb F_p((t))$ 
one has $v_F(\alpha^{\sigma}-\alpha)\ge v_F(\alpha)+i(\sigma)$
with equality when 
$v_F(\alpha)$ is relatively prime to  $p$. Therefore
$$i(\sigma^p)=v_F((\sigma^p-1)(t))-1=v_F((\sigma-1)^p(t))-1
\ge pi(\sigma),$$ 
hence $R_n^p\le R_{np}$  and 
$i(\sigma^p)=v_F((\sigma-1)^p(t))-1=pi(\sigma)$ for $i(\sigma)=1+qi$.

(2) The proof of this part is rather lengthy.

\proclaim{\bdas Lemma 1}\idas
Let  $j=\overline{\jmath}p^n$, $i=\overline{\imath}p^m$ with integers $\overline{\imath},\overline{\jmath}$ relatively prime to
$p$. 

\noindent Let $o\in\Bbb N$ be such that $i+j, 2i$ are equal or greater than $o$.
Then for $a,b\in\Bbb F_p$ 
$$
\align
&[t+at^{1+qj},t+bt^{1+qi}]
\equiv t+\overline{\jmath}abt^{1+qj+q^2ip^n}/(1+at^{qj})\\
&-bt^{1+qi}
\bigl(\sum_{l=1}^{\overline{\imath}}\binom{\overline{\imath}}
la^l(t^{q^2ljp^m}/(1+at^{q^3j})^{ljp^m})\bigr)/\bigl((1+bt^{qi})(1+at^{q^2j})^i
\bigr) \mod t^{1+q^2o+q}\Bbb F_p[[t]].
\endalign
$$
\endproclaim
\demo {\idas Proof}\das
Let 
$g(t)=(t+at^{1+qj})\circ(t+bt^{1+qi})-(t+bt^{1+qi})\circ(t+at^{1+qj})$.
Then 
$$%\split
g(t)=at^{1+qj}(1+bt^{qi})\biggl(\sum_{l=1}^{\overline{\jmath}}\binom{\overline{\jmath}}l b^lt^{lq^2ip^n}\biggr)%\\
-bt^{1+qi}(1+at^{qj})\biggl(\sum_{l=1}^{\overline{\imath}}\binom{\overline{\imath}}
la^lt^{lq^2jp^m}\biggr).
%\endsplit
$$
Denote by $\alpha(t)$ and $\beta(t)$ the inverse series in $T$ to
$t+at^{1+qj}$ and $t+bt^{1+qi}$ respectively.
Note that for $f_1(t)\in 1+t^q\,\Bbb F_p[[t^q]]$, $f_2(t)\in
1+t^l\,\Bbb F_p[[t]]$
$$(tf_1(t))\circ(tf_2(t))\equiv tf_1(t)f_2(t)\mod t^{1+q(l+1)}\,\Bbb F_p[[t]].
\tag"$(\triangledown)$"$$

Therefore,
$$\aligned
&\alpha(t)\equiv t/(1+at^{qj})\mod t^{q^2j}\,\Bbb F_p[[t]],\qquad
\beta(t)\equiv t/(1+bt^{qi})\mod t^{q^2i}\,\Bbb F_p[[t]],\\
&(\beta\circ\alpha)(t)\equiv \alpha(t)/(1+bt^{qi})\equiv t/((1+bt^{qi})
(1+at^{qj}))\mod t^{q^2j}\,\Bbb F_p[[t]],\\
&(\alpha\circ\beta)(t)\equiv \alpha(t)/(1+bt^{qi}) \mod t^{q^2i}\,\Bbb F_p[[t]].
\endaligned
$$
Due to the conditions on $i,j,o$
the series $g\circ\alpha\circ\beta$ is congruent $\mod t^{1+q^2o+q}\,\Bbb F_p[[t]]$ to
$$\aligned &\bigl(at^{qj}\overline{\jmath}bt^{q^2ip^n}(t+bt^{1+qi})\bigr)\circ
\bigl(\alpha(t)/(1+bt^{qi})\bigr)\\
&-\bigl(bt^{qi}\bigl(\sum_{l=1}^{\overline{\imath}}\binom{\overline{\imath}}
la^lt^{lq^2jp^m}\bigr)(t+at^{1+qj})\bigr)\circ
\bigl(\alpha(t)/(1+bt^{qi})\bigr)\\
&\equiv 
\overline{\jmath}abt^{qj+q^2ip^n}\alpha(t)
-b\alpha(t)^{qi}
\bigl(\sum_{l=1}^{\overline{\imath}}\binom{\overline{\imath}}la^l(t/
(1+at^{1+qj}))^{lq^2jp^m}\bigr)
\beta(t)
\endaligned
$$
and the lemma follows.
\qed\enddemo 

\proclaim{\bdas Remark}
\idas
If $i\ge j$, one can take $i+j$ as $o$. 
If a pair $(i,j)$ satisfies the conditions of Lemma 1,
then so does every pair $(i',j')$ with $i'\ge i, i'+j'\ge i+j$.
\endproclaim

\das 
Now it is clear that $[t+at^{1+qj},t+bt^{1+qi}]$ for $i> j$
belongs to $t+t^{1+q(i+qj)}\,\Bbb F_p[[t]]$
and $[T_i,T_i]\le T_{(q+1)i+1}$.
Then (1) implies that
$T_i/T_{(q+1)i}$ is abelian of exponent $pq$.

\proclaim{\bdas Lemma 2}\idas
Let $p>2$. Fix $s$ satisfying $1\le s \le r$. 
Let $i>j\ge q^2$ and $i$ be relatively prime to $p$.  
Let the set $X$ consist
of pairs $(i_m,j_m)$ of the following type{\das:}

\noindent $i_m\ge i$, $j_m\ge j-q$;

\noindent $i_m$ is relatively prime to $p$;

\noindent $j_m\ge j$ if $i_m=i$; 
$qj_m+p^si_m>qj+p^si$ if $i_m>i$;

\noindent  $(\diamond)\qquad \text{\idas if $i_m+qj_m<j+qi$, then $i_m=u_mi+v_mq$ with integers }u_m\ge 1, v_m\ge 0.$

Put  $j_m=\overline{\jmath}_mp^{n(j_m)}$ where $\overline{\jmath}_m$
is relatively prime to $p$.
Let $c_m,d_m,e_m,f_m$ be non-negative integers such that
$d_m>0$ iff $c_m>0$ and $f_m>0$ iff $e_m>0$.

Then the equality
$$\sum_{(i_m,j_m)\in X} ( c_m{i_m}+qd_mj_m+e_m\overline\jmath_m+qf_mi_mp^{n(j_m)})
=k+qj, \quad p^{s-1}i<k\le p^si, p^s|k $$ 
implies that up to reodering terms 
{\das (1)} $k=p^si${\das; (2)} if $ p^s< q$ then
$c_1=p^s,d_1=1$, $i_1=i,j_1=j$, and all $e_m,f_m$ are 
zero{\das; (3)} if $ p^s=q$, then either $c_1=q,d_1=1$, $i_1=i,j_1=j$, and all $e_m,f_m$ are zero
or $e_1=q, f_1=1$, $i_1=i, j_1=j$ and all $c_m,d_m$ are zero.
\endproclaim
\demo {\idas Proof}\das
Let not all $c_m$ be zero. Since $i_m+qj_m>qj$,
we have $\sum (c_m{i_m}+qd_mj_m)=k+qj$.
If, say, $c_1i_1+qd_1j_1\ge j+qi$, then
$c_1i_1+qd_1j_1= k+qj$, $p^s|c_1$ and
$c_1i_1+qd_1j_1\ge p^si_1+qj_1$ (which is $>p^si+qj$ if $i_1\not=i$).
Hence $i_1=i,j_1=j, c_1=p^s, d_1=1, k=p^si$.
If all $c_mi_m+qd_mj_m$ are smaller than $j+qi$, then
$\sum c_mi_m=(\sum c_mu_m)i+(\sum c_mv_m)q$,
$p^s$ divides $\sum c_mu_m$ and $\sum c_mi_m\ge p^si$ by $(\diamond)$.
Then $k+qj\ge p^si+q(\sum d_mj_m)$.
Since $2j_m>j$, we deduce that, say, $d_1=1$, $d_m=0$ for $m\not=1$.
Then $c_1i_1+qj_1=k+qj$.
Since  $p^s$ divides $c_1$, $c_1i_1+qj_1\ge p^si_1+qj_1$, and the latter is
$>p^si+qj$ whenever $i_1\not=i$. Hence $i_1=i,j_1=j, c_1=p^s, k=p^si$.

Let all $c_m$ be zero. Then, say, $f_1=1$ and $f_m=0$ for $m\not=1$. 
From $e_1\overline{\jmath}_1+qi_1p^{n(j_1)}=k+qj$
it follows that $n(j_1)=0$ and $p^s$ divides $e_1$.
Then $k+qj\ge p^sj_1+qi_1$, therefore $j_1<i_1$. In addition
$p^sj_1+qi_1> p^si_1+qj_1$ if $p^s<q$,
and $p^si_1+qj_1>p^si+qj$ if 
$i_1\not=i$. Therefore,  
$i_1=i, j_1=j, e_1=p^s=q, k=qi$.
\qed\enddemo 

\das 

From the previous Lemmas and formula $(\circ)$ of sect. 1 
we deduce the following

\proclaim{\bdas Remark}
\idas Let $i,j$ be as in Lemma 2.
For $i$  relatively prime to $p$ and $1\le s <r$ the coefficient of $t^{1+q(p^si+qj)}$ in 
$[t+at^{1+qj}+\dots,t+bt^{1+qi}]$ is  equal to $-iab$
and the coefficient of $t^{1+q^2(i+j)}$  is equal to
$(j-i)ab$.
\endproclaim

\das 

Now let $H$ be a nontrivial closed normal subgroup of an open subgroup $G$ of $T$. 
We say that a natural number $l$ belongs to $H$ if
$T_l\le T_{l+1}H$.
From what has been proved it follows that
given $n\ge0$  there is $l_n$ such that all $p^nl\ge l_n$ with $l$ relatively prime to $p$ belong to $H$. 

Let $j\ge \max(l_0+q,q^2)$, $i>j$,  and let $i,i-j$ be relatively prime to $p$. 
Let $T_i\le G$.

We shall use Lemma 1, Lemma 2, two previous Remarks and formula $(\circ)$ of sect. 1 in the following arguments.

Let $\theta(t)=t+at^{1+qj}+\dots\in H$.
Then 
$$\omega_0(i,j,i+qj)=[\theta(t),t+t^{1+qi}]\equiv
t+\sum_{i_m\ge i}w_{i_m}t^{1+q(qj+i_m)}\mod t^{1+q^2(i+j)+q}\,\Bbb F_p[[t]]$$
where $(i_m,j)$ satisfy $(\diamond)$ and
$w_{p^si}\not=0$ for  $1\le s \le r$.
Let $\rho(t)=t+at^{1+q(j-p)}+\dots\in H$.  Then
using $(\triangledown)$ we deduce that 
for appropriate $c\in\Bbb F_p$
$$%\split
[\theta(t),t+t^{1+qi}][\rho(t),t+ct^{1+q(i+qp)}]%\\
\equiv
t+\sum_{i_m>i} x_{i_m}t^{1+q(qj+i_m)}\mod t^{1+q^2(i+j)+q}\,\Bbb F_p[[t]]
%\endsplit
$$
where $(i_m,j)$ satisfy $(\diamond)$ and
$x_{p^si}\not=0$ for  $1\le s \le r$.

To deal with $t+t^{1+q(qj+i_m)}$ where $p^{s-1}i<i_m< p^si$,
$i_m=p^n\overline{\imath}_m$ and $\overline{\imath}_m>i$ is relatively prime to $p$ note that by Lemma 2 and induction $0\le n\le s-1$.
Then use appropriate 
$\omega_n(\overline{\imath}_m,j,i_m+qj)$ for $i_m+qj\ge j+qi$
and $\omega_n(u_mp^{-n}i,j+v_m,i_m+qj)$ for $i_m+qj<j+qi$, 
$i_m=u_mi+v_mq$.
Hence in $H$ there is an element 
$$\omega_s(i,j,p^si+qj)\equiv
t+\sum_{i_m\ge p^si} y_{i_m}t^{1+q(qj+i_m)}\mod t^{1+q^2(i+j)+q}\,\Bbb F_p[[t]]$$
where $(i_m,j)$ satisfy $(\diamond)$ and
$y_{p^li}\not=0$ for $s\le l\le r$.
To eliminate $t+t^{1+ q(qj+p^si)}$ if $s<r$ multiply 
$\omega_s(i,j,p^si+qj)$ by appropriate power of $\omega_s(i+q,j-p^s,p^si+qj)$.

We conclude that the numbers
$p^si+q^2j$ for $1\le s\le r$ belong to $H$.
Varying $i$ for the given $j$,
then passing to a greater $j'$ such that the pair
$(i',j')$  satisfies the conditions
above and varying $i'$
we conclude  that all sufficiently large $l$ belong to $H$.
Since $H$ is closed,
$T_l$ is contained in $H$ for sufficiently large $l$
and $H$ is open.

(3) Note that if $i>j$ and $i$ is relatively prime to $p$,
then $[t+at^{1+qj},t+bt^{1+qi}]=t-iabt^{1+qi+q^2j}+\dots$.
Therefore $T_{l}\le T_{l+1}[T,T]$ for all $l\ge q+2$
relatively prime to $p$.
By (1) $T_{pi}\le T_{pi+1}T^p$ and the assertion follows.

(4) The group $T$ doesn't contain an open subgroup of finite rank,  
since the cardinality of 
%!
$T_i/[T_i,T_i]T_i^p$ tends to infinity when $i$ tends
to infinity.  Alternatively, 
if $T$ were a $p$-adic analytic group,
then by [6, Th. 9.34]
it would contain 
an open subgroup $G$ which is
a uniformly powerful pro-$p$-group,
so in particular $G^p=\{g^p:g\in G\}$ would be a subgroup of $[G,G]$,
see [6, Prop. 2.6]. However,  the subgroup
$[G,G]$ contains some elements $t+t^{1+qi}+\dots$ with $i$ relatively prime to $p$ which are obviously not in $G^p$.\qed\enddemo 
\das

\proclaim{\bdas Remark 1} \das
Using the commutator formula of Lemma 1 one can show that $T$ has finite width,
i.e. the orders of terms of the quotient filtation of the
lower central series filtration of $T$ are uniformly bounded (see [18]).
\endproclaim
%\comment
\demo {\idas Proof}\das 
Let $\gamma_i(T)=[T,\gamma_{i-1}(T)]$ and $\gamma_1(T)=T$.
From the formula of Lemma 1 it follows that
(1) the minimal natural number which belongs to $\gamma_i(T)$ is
$2+(i-1)q$; all greater natural numbers relatively prime to $p$ belong to
$\gamma_i(T)$;
(2) if $l$ belongs to $\gamma_{i+1}(T)$ and $l$ is divisible by $q$,
then $l\ge q(3+(i-1)q)$; all  $l\ge (5+(i-1)q)q$ divisible by $q$ belong to 
$\gamma_{i+1}(T)$;
(3) if $l$ belongs to $\gamma_{i}(T)$ and $l$ is not divisible by $q$,
then $l+q^2$ belongs to  $\gamma_{i+1}(T)$.

Therefore the width of $T$ which is equal to $\sup\log_p|\gamma_{i}(T):\gamma_{i+1}(T)|$ is not greater than $q^2+2$.
\enddemo
%\endcomment

\proclaim{\bdas Remark 2} \das
Let $(A,\Cal M)$ be a commutative complete local noetherian
integral domain with finite residue field of characteristic $p$.
For the theory of analytic groups over $A$ see [21], [27, Part IV].
Using the commutator formulas of Lemma 1 one can show that the wild group
$R$ is not analytic over $A$, see [10].
\endproclaim

\das
\subhead 7 \endsubhead 
Now let $q=p^r>p$ and $F$ be an unramified extension 
of $\Bbb Q_p$ of degree $\ge q$.
Let $L/F$ be a Galois totally ramified extension with the Galois group
isomorphic to $T$ (it exists, since according
to Shafarevich's theorem the pro-$p$-part of the absolute
Galois group of $F$ is a free pro-$p$-group with $|F:\Bbb Q_p|+1$ generators).

\smallskip

A Galois extension $L/F$
is called deeply ramified if the set of its upper ramification
breaks is unbounded
(see [5] and some additional information in [8]).

\proclaim{\bdas Theorem} \idas 

{\das (1)} The extension $L/F$ is deeply ramified.

{\das (2)} The extension $L/F$ is arithmetically profinite{\das;} 
if $F\subset K\subset E\subset L$
with finite $K/F$ and infinite Galois $E/K$,
then $G(E/K)$ is not  $p$-adic analytic.
\endproclaim \das
\demo {\idas Proof}\das 
(1) (inspired by Sen's proof [24] of Serre's
conjecture on the ramification filtration of $p$-adic analytic extensions
of local number fields)

Denote by $F_i$ the fixed field of $T_i$.
For a finite Galois extension $E/K$ and an automorphism $\sigma\not=1$ of $E/K$ denote
by $t_K(\sigma)$ the maximal rational number $x$ such that
$\sigma$ doesn't belong to  the upper ramification group
$G(E/K)^v$ for every $v>x$.
Denote  by $u(E/K)=\max \{t_K(\sigma): \sigma\in G(E/K),
\sigma\not=1\}$ its maximal upper
ramification break
(see [9, Ch. III, sect. 5]).
Due to local class field theory for a finite totally ramified abelian $p$-extension  $M/K$ the number $u(M/K)$ 
is equal to the maximal integer $n$ such that 
$U_{n+1,K}N_{M/K}U_{1,M}$ doesn't contain
the $n$th group of principal units $U_{n,K}$ of $K$.

We shall work with abelian extension inside $L/F$.
Note that $L$ doesn't contain 
nontrivial $p$th roots of unity.
Following Sen [24] we say that the abelian extension $M/K$ is not small
if $u(M/K)\ge pe(K)/(p-1)$ where $e(K)$ is the absolute ramification index
of $K$. Then, due to local class field theory there is an automorphism
$\sigma\in G(M/K)$ such that $t_K(\sigma)\ge e(K)/(p-1)$ and for every 
automorphism with this property we get
$t_K(\sigma^p)=t_K(\sigma)+e(K)$.

The extension $F_{(q+1)i}/F_i$ is abelian.
Denote $i'=\lceil (k(q+1)/p\rceil $. Then due to property (1) of 
the proposition the group
$G(F_{(q+1)i}/F_{i'})$ has exponent $p$ and each automorphism of
$F_{(q+1)i}/F_i$ which is nontrivial on $F_{i'}$
is of exponent at least $p^2$.

Suppose that $F_{(q+1)i}/F_i$ is not small.
Then there is $\sigma$ in $G(F_{(q+1)i}/F_i)$ satisfying the property 
$$t_{F_i}(\sigma^p)=t_{F_i}(\sigma)+e(F_i).\tag"\das(*)"$$
Note that if $\sigma^{p^2}\not=1$ in $G(F_{(q+1)i}/F_i)$,
then $\sigma^p$ satisfies (*);
if $1\not=\tau\in G(F_{(q+1)i}/F_i)$ with
$t_{F_i}(\tau)>t_{F_i}(\sigma^p)$, then
$t_{F_i}(\tau)>pe(F_i)/(p-1)$ and hence $\tau=\rho^p$ for some
$\rho\in G(F_{(q+1)i}/F_i)$ and $\rho$ satisfies (*).

Let $1\not=\sigma\in G(F_{(q+1)i}/F_i)$ be such that
$u(F_{(q+1)i}/F_i)=t_{F_i}(\sigma^p)=t_{F_i}(\sigma)+e(F_i)$
and such that $\sigma$ acts trivially on $F_{k}$
with the maximal possible integer $k=k(i)$.
We get $i'>k\ge i'/p>i$.
Denote $\sigma_k=\sigma$.
Since $t_{F_i}(\sigma_k)\ge t_{F_i}(\tau)$ for every $\tau$
which acts nontrivially on $F_{i'}$, we deduce that
$t_{F_i}(\sigma_k)=u(F_{k+1}/F_i)$ and therefore 
$t_{F_k}(\sigma_k^p)=t_{F_k}(\sigma_k)+e(F_k)$.

Since $t_{F_k}(\sigma_k^p)>e(F_k)/(p-1)$,
and $\sigma_k^{p^2}$ is nontrivial in $G(F_{(q+1)k}/F_k)$,
we deduce that $t_{F_k}(\sigma_k^{p^2})=t_{F_k}(\sigma_k^p)+e(F_k)$
and the extension $F_{(q+1)k}/F_k$ is not small.

Start with the extension $F_{q+1}/F_1$ which is abelian 
of exponent $p^3$. 
Since $e(F_1)=1$, raising to the $p$th
power maps $U_{i,F_1}$ onto $U_{i+1,F_1}$, therefore
$u(F_{q+1}/F_1)\ge 3$ and $F_{q+1}/F_1$ is not small.
Put $k_1=k(1)\ge q/p^2+1$, $k_{n+1}=k(k_{n})>k_{n}$.
Then $t_{F}(\sigma_{k_{n+1}})\ge
t_{F}(\sigma_{k_n}^p)\ge t_{F}(\sigma_{k_n})+e(F)$.
Thus, $q(F_{k_n}/F)$ tends to $+\infty$ when $n$ tends to $+\infty$,
hence $L/F$ is deeply ramified.

(2) The extension $L/F$ is arithmetically profinite:
every $G(L/F)^x$ is a nontrivial normal closed
subgroup of $T$, therefore is open of finite index by property (2)
of the proposition.
Normal extensions $E$ of $K$ in $L$ are either finite over $K$
or coincide with $L$ and therefore $G(E/K)$ is not
$p$-adic analytic by the previous proposition. 
\qed\enddemo 
\das

J. Coates and R. Greenberg  in 
[5, p. 144] stated the following problem:
is it true that for every finite extension $K$ of $\Bbb Q_p$
there exists a deeply ramified Galois $p$-extension $M$ of $K$
such that that no subfield $M'$ of $M$ is an infinite Galois extension
of a finite extension $Q$ of $\Bbb Q_p$ with
$G(M'/Q)$ being $p$-adic analytic.
%This problem is related to Fontaine--Mazur's conjecture which says
%that there do not exist a number field which has
%an infinite unramified everywhere Galois $p$-extension
%which has infinite $p$-adic analytic quotients.

Using the previous theorem one can provide the affirmative answer on 
Coates--Greenberg's problem by  taking 
the normal closure of $KL$ over $K$ as $M$.
Indeed, then $M/K$ is deeply ramified.
Note that $M$ is the compositum of a finite number of fields $\sigma (KL)$
and denote the compositum of all of them with $KL$ excluded by $N$. 
If $KL\not\subset N$, then $KL\cap N$ is of finite degree over
$KF$, since $G(KL/KF)$ is hereditarily just infinite.
If there is a subfield $M'$ of $M$ which is a
Galois extension of a finite extension $Q$ of $\Bbb Q_p$
with $G(M'/Q)$ being $p$-adic analytic, then,
since $QKL\cap M'$ is of finite degree over $QKF$ ($G(QKL/QKF)$ is hereditarily just infinite), the group $G(M'KL/QKL)$ isomorphic to $G(M'/M'\cap QKL)$ 
is an infinite $p$-adic analytic group, 
so $G(LN/L)$ has an open subgroup
which has  an infinite $p$-adic analytic quotient. 
Since $N\cap L$ is of finite degree over $F$, $G(N/F)$ has an open
subgroup which has an infinite $p$-adic analytic quotient. Eventually we deduce that $G(L/F)$  has an open
subgroup which has an infinite $p$-adic analytic quotient, a contradiction. 

Thus, $G(M/K)$ doesn't have 
open subgroups which have infinite $p$-adic analytic quotients
and $M$ is an arithmetically profinite $p$-extension of $K$.

The affirmative answer on the problem stated in sect. 4 implies
that one can realize the Nottingham group as the Galois group
of an arithmetically profinite extension of a local number field,
therefore providing another extension $L/F$ which has all the properties
of the theorem.

\das
\subhead 8 \endsubhead 
Let $(A,\Cal M)$ be a commutative complete local noetherian
integral domain with finite residue field of characteristic $p$.
Let $\nu (\alpha)=\min \{i: \alpha\in \Cal M^i\}$.

Let $F=F(X,Y)=X+_F Y$ be a formal noncommutative group law over $A$ of dimension $d$
and $F((\Cal M)^d)$ the group of formal points associated with $F$.
Let
$$X+_F Y+_F-_F X-_F Y=\sum_{|J|+|K|\ge 1}
(c_{J,K}^{(i)}X^JY^K)_{1\le i\le d}$$
with the usual notations for multi-indices.

Let $m=\min\{|J|+|K|:c_{J,K}^{(i)}\not=0\text{\das \  \ for some $i$}\}$.

\noindent Let $j'=\max\{|J|:c_{J,K}^{(i)}\not=0 \text{\das \ \ for some $i$}, |J|+|K|=m\}$.
\noindent Let $(J',K')$ be the maximal with respect to the lexicographical order
(the first component is the most respected) multi-index such that
$|J'|=j', |K'|=m-j'$ and $c_{J',K'}^{(i)}\not=0$ for some $1\le i\le d$.
Let $i_0$ be the minimal index such that $c_{J',K'}^{(i_0)}\not=0$.
Denote $s=\nu(c_{J',K'}^{(i_0)})$.

Denote by $B_n$ the set of vectors $(b_1,\dots,b_d)$ of length $d$
satisfying $\nu(b_i)\ge n+i(s+2)$.
Then  $F(B_n)$ is an open subgroup of $F((\Cal M)^d)$ for $n>(d-2)(s+2)$.

\proclaim{\das Lemma}\idas Let $c>d(s+2)$ and $n>m(c+d(s+2))$.
Then 
for every $\alpha\in F(B_n)\setminus F(B_{n+1})$, $\beta\in F(B_{n+c})\setminus
F(B_{n+c+1})$
$$[\alpha,\beta]=\alpha+_F\beta-_F\alpha-_F\beta\not\in F(B_{2+s+m(c+d(s+2))+mn}).$$
\endproclaim

\demo{\idas Proof}\das
Indeed, $\nu(c_{J',K'}^{(i_0)}\alpha^{J'}\beta^{K'})
= s+J'A+K'B\le s+m(n+c+d(s+2)+1)$
where $A=(\nu(a_1),\dots,\nu(a_d))$ for $\alpha=(a_1,\dots,a_d)$, $B$ for $\beta$ is defined similarly.

If $|J|=j', |K|=m-j'$ and $(J,K)<(J',K')$ then $$\nu(c_{J,K}^{(i_0)}\alpha^{J}\beta^{K})
\ge s+1+J'A+K'B>\nu(c_{J',K'}^{(i_0)}\alpha^{J'}\beta^{K'}).$$
If $m=|J|+|K|$ and $|J|<j'$, then
$$\nu(c_{J,K}^{(i_0)}\alpha^{J}\beta^{K})\ge c-(d-1)(s+2)-1+J'A+K'B>
\nu(c_{J',K'}^{(i_0)}\alpha^{J'}\beta^{K'}).$$
If $m<|J|+|K|$, then for appropriate $J'',K''$ with $|J''|+|K''|=m$
$$\split\nu(c_{J,K}^{(i_0)}\alpha^{J}\beta^{K})\ge 
n+s+2+J''A+K''B
\ge (n+s+2)(1+m)\\>s+m(n+c+d(s+2)+1)\ge s+J'A+K'B.
\endsplit$$
The lemma is proved.
\enddemo

\das 

\proclaim{\das Proposition}\idas The wild group $R$ is not $A$-analytic.
\endproclaim
\demo{\idas Proof} \das Assume that  $R$ is $A$-analytic and get a contradiction.
Let $F$ be the corresponding formal group of dimension $d$ and
let $R$ contain $F(B_o)$ as an open subgroup, see [7, Th. 22].
Since $R_i/R_{2i}$ is isomorphic to $(\Bbb F_p)^i$,
one can choose some of $R_j$ for a new filtration $F(B_o)\ge W_1\triangleright W_2\triangleright \dots$
such that the quotient groups $W_i/W_{i+1}$ are abelian of the same cardinality as $(A/\Cal M)^d$. 

The wild group is finitely generated, therefore 
the topology  $\{F(B_i)\}$ is equivalent to the topology $W_j$
and so there is $z$ such that $W_1\ge F(B_z)$. Then every $F(B_i)$ for
$i\ge z$ is normal in $W_1$.

It is easy to check that if $W_1\ge R_v$ then $(R_j\setminus R_{j+1})\cap H\not=\varnothing\Rightarrow H\ge R_{j+v+2p}$.
From here one deduces that there is $u\in\Bbb N$ such that for every normal subgroup $H$ of $W_1$
$$H\cap (W_j\setminus W_{i+1})\not=\varnothing\Rightarrow H\ge W_{j+u}.$$

Let $j_1(i)=\min \{j: (W_j\setminus W_{j+1})\cap(F(B_i)\setminus F(B_{i+1}))
\not=\varnothing\}$. Then 
$$W_{j_1(i)+u}\le F(B_i)\qquad\text{\das for $i\gg0$}.\tag 1$$
The sequence $j_1(i)$ is non-decreasing.

Let $j_2(i)=\max \{j: (W_j\setminus W_{j+1})\cap(F(B_i)\setminus F(B_{i+1}))
\not=\varnothing\}$. Then 
$$j_1(i)\le j_2(i)\le j_1(i+1)+u\qquad\text{\das for $i\gg 0$.}\tag 2$$
If $j_1(i)>1+(i-o)$, then $F(B_i)\le W_{i-o+2}$ and the index
of $F(B_i)$ in $F(B_o)$ is smaller than that of $W_{i+o-2}$ in $W_1$,
a contradiction. Hence by (1) 
$$j_1(i)\le 1+i-o\quad\text{\das  and }\quad W_{1+i-o+u}\le F(B_i) \text{\das
\ \  for $i\gg0$}.\tag 3$$

If $W_j\le F(B_i)$, then $j-1\ge i-z$, hence
$W_{i-z}\cap (F(B_o)\setminus F(B_i))\not=\varnothing$.
Let $\alpha\in (F(B_{i'})\setminus F(B_{i'+1}))\cap W_{i-z}$ 
for some $i'<i$. Then $j_2(i')\ge i-z$ and by (2) $j_1(i)+u\ge j_1(i'+1)+u
\ge j_2(i')$, hence 
$j_1(i)\ge i-z-u$.

Thus,
$$i-z-u\le j_1(i)\le j_2(i)\le j_1(i+1)+u\le 2+i-o+u\qquad\text{\das for $i\gg 0$}.\tag 4$$

Let $\alpha\in (W_j\setminus W_{j+1})\cap(F(B_n)\setminus F(B_{n+1}))$ and
$\beta\in (W_{j+g}\setminus W_{j+g+1})\cap(F(B_{n+c})\setminus F(B_{n+c+1}))$.
Then by (4) 
$$j+o-u-2\le n\le j+z+u, \quad j+g+o-u-2\le n+c\le j+g+z+u.\tag 5$$
Hence
$$g-2u+o-z-2\le c\le g+2u+z-o+2.\tag 6$$

Choose $q=p^r$ and $i$ sufficiently large such that
$$j\ge 3-o+u+s+m(g+3u+2z-o+d(s+2)+2), g-2u+o-z> d(s+2)+2$$
and 
$\alpha=t+t^{1+qi}\in W_j\setminus W_{j+1}$, $\beta=t+t^{1+qi+q}\in W_{j+g}\setminus W_{j+g+1}$,
$R_{1+q^2i+q(i+1)}\le W_{(m+1)j}$.
Then $c>d(s+2)$ by (6) 
and $n\ge j+o-u-2>m(g+2u+z-o+2+d(s+2))\ge m(c+d(s+2))$,
$(m+1)j-1+o-u\ge 2+s+m(g+2u+z-o+2+d(s+2)+j+z+u)\ge
2+s+m(c+d(s+2)+n)$ by (5).
According to sect. 6 
we get $[t+t^{1+qi},t+t^{1+q(i+1)}]\in R_{1+q^2i+q(i+1)}\le W_{(m+1)j}$
and the latter is a subgroup of $F(B_{(m+1)j-1+o-u})
\le F(B_{2+s+m(c+d(s+2)+n)})$ by (3).

Thus, there are 
$$\alpha\in W_n\setminus W_{n+1}, \beta\in W_{n+c}\setminus W_{n+c+1},
\quad\text{\das such that}\quad
[\alpha,\beta]\in F(B_{2+s+m(c+d(s+2)+n)}).$$ 
That contradicts the property
of $A$-analytic groups in the lemma. 
\enddemo

\vskip .6in

\tenpoint

\centerline{{\bdas  References}}
\das
\vskip .5in

%\iit{#}{\item"{#.}"}
\roster

\item"{[1]}" Barnea Y., Shalev, A.:
Hausdorff dimension, pro-$p$-groups, and Kac--Moody algebras.
To appear in Trans. of the AMS, 1997

\item"{[2]}" Barnea Y.:
The spectrum of the Nottingham group.
Preprint, Jerusalem 1997

%\vfill\pagebreak

\item"{[3]}" Camina, R.:
Subgroups of the Nottingham group.
J. Algebra {\bf 196}, 101--113(1997)

\item"{[4]}" Cartwright, D.:
Pro-$p$-groups and Lie algebras over local fields of characteristic $p$.
Ph. D. Thesis, QMW, London 1997

\item"{[5]}" Coates, J., Greenberg, R.:
Kummer theory for abelian varieties over local fields.
Invent. Math. {\bf 124}, 129--174(1996) 

\item"{[6]}" Dixon, J.D., du Sautoy, M.P.F.,Mann, A., Segal, D.:
Analytic pro-$p$-groups.
LMS Lect. Notes Series {\bf 157}
Cambridge Univ. Press, Cambridge 1991

\item"{[7]}"  du Sautoy, M.P.F.:
Pro-$p$-Groups.
Preprint, Cambridge 1997

\item"{[8]}" Fesenko, I.:
On deeply ramified extensions.
Preprint, Heidelberg 1995

\item"{[9]}" Fesenko, I., Vostokov, S.:
Local Fields and Their Extensions.
AMS, Providence, R.I. 1993

\item"{[10]}" Fesenko, I. B.:
The wild group is not analytic over complete local noetherian rings
with finite residue field.
http://www.maths.nott.ac.uk/personal/ibf/non.ps,
Nottingham 1997

\item"{[11]}" Fontaine, J.--M.:
Un r\'esultat de Sen sur les automorphismes de corps locaux.
S\'eminaire Delange-Pisot-Poitou, Paris (1969--1970) expos\'e 6

\item"{[12]}" Iwasawa, K.:
On solvable extensions of algebraic number fields.
Ann. Math. {\bf 58}, 548--572(1953)

\item"{[13]}" Johnson, D.L.:
The group of formal power series under substitution.
J. Austral. Math. Soc. (Series A) {\bf 45}, 298--302(1988)

\item"{[14]}"  Klopsch, B.:
Normal subgroups and automorphisms of substitution groups 
of formal power series.
Preprint, Oxford 1997

\item"{[15]}" Laubie, F.:
Sur la ramification des extensions de Lie.
Comp. Math. {\bf 98}, 253--262(1985)

\item"{[16]}" Laubie, F.:
Extensions de Lie et groupes d'automorphismes de corps locaux.
Comp. Math. {\bf 67}, 165--189(1988)

\item"{[17]}"  Laubie, F., Saine, M.:
Ramification of automorphisms of $k((t))$.
J. Number Theory {\bf 63}, 143--145(1997)

\item"{[18]}" Leedham--Green, C.R., Plesken, W., Klaas, G.:
Linear pro-$p$-groups of finite width.
Preprint,  Aachen--London 1997

\item"{[19]}" Li, H.-Ch.:
$p$-adic periodic points and Sen's theorem.
J. Number Theory {\bf 56}, 309--318(1996)

\item"{[20]}" Lubin, J.:
Non-archimedean dynamical systems.
Comp. Math. {\bf 94}, 321--346(1994)

\item"{[21]}" Lubotzky, A., Shalev, A.:
On some $\Lambda$-analytic pro-$p$-groups.
Israel J. Math. {\bf 85}, 307--337(1994)

\item"{[22]}" Lubotzky, A., Wilson, J.S.:
An embedding theorem for profinite groups.
Arch. Math. {\bf 42}, 397--399(1984)

\item"{[23]}" Sen, Sh.:
On automorphisms of local fields.
Ann. Math. {\bf 90}, 33--46(1969)

\item"{[24]}" Sen, Sh.:
Ramification in $p$-adic Lie extensions.
Invent. Math. {\bf 17}, 44--50(1972)

\item"{[25]}" Wintenberger, J.--P.:
Automorphismes et extensions galoisiennes de corps locaux.
Th\`ese de 3 cycle, Grenoble 1978

\item"{[26]}" Wintenberger, J.--P.:
Extensions de Lie et groupes d'automorphismes
des corps locaux de caract\'eristique $p$.
C. R. Acad. Sc. {\bf 288} s\'erie A, 477-479(1979)

\item"{[27]}" Wintenberger, J.--P.:
Extensions ab\'eliennes et groupes d'automorphismes des corps locaux.
C. R. Acad. Sc. {\bf 290} s\'erie A, 201--203(1980)

\item"{[28]}" Wintenberger, J.--P.:
Le corps des normes de certaines extensions infinies des corps locaux; applications.
Ann. Sci. E.N.S., 4 s\'erie {\bf 16}, 59--89(1983)

\item"{[29]}" York, I. O.:
The group of formal power series under substitution.
Ph. D. Thesis, Nottingham 1990

\bigskip
\bigskip

Department of Mathematics

University of Nottingham

NG7 2RD Nottingham England

\bye